\def \verbose {0}

\documentclass[12pt]{article}
\input epsf      
\usepackage{graphicx}
\usepackage{amssymb}
\usepackage{amsfonts}
\usepackage{amsthm}
\usepackage{amsfonts}

\newcommand \slope{{0.520655}}

\newcommand \half {{\mathbb H}}
\newcommand \disc {{\mathbb D}}

\newcommand{\integers}{{\mathbb Z}}

\newcommand{\compose}{{\circ}}

\newcommand \no {\noindent}
\newcommand \ra {\rightarrow}

\newcommand{\ba}[1]{\begin{array}{#1}}
\newcommand{\ea}{\end{array}}

\newcommand{\be}{\begin{equation}}
\newcommand{\ee}{\end{equation}}
\newcommand{\bea}{\begin{eqnarray}}
\newcommand{\eea}{\end{eqnarray}}
\newcommand{\beann}{\begin{eqnarray*}}
\newcommand{\eeann}{\end{eqnarray*}}

\newcommand{\sech}{{\rm sech }}

\newcommand{\bp}{{\overline{b}}}

\newcommand{\prob}{{\bf P}}

\newcommand{\probNdelta}{{\bf P}_{\half,N,\delta}}
\newcommand{\probfl}{{\bf P}_\half^{fixed}}
\newcommand{\probzdelta}{{\bf P}_{\half,z,\delta}}
\newcommand{\probNp}{{\bf P}_{N}^\prime}

\newcommand{\probsle}{{\bf P}^{radial}_{\half,i}}

\newcommand{\EN}{{\bf E}_N}
\newcommand{\ENdelta}{{\bf E}_{\half,N,\delta}}
\newcommand{\Efl}{{\bf E}_\half^{fixed}}

\def\reff#1{(\ref{#1})}
\begin{document}

\bibstyle{ams}

\title{Transforming fixed-length self-avoiding walks \\
into radial SLE$_{8/3}$}

\author{Tom Kennedy
\\Department of Mathematics
\\University of Arizona
\\Tucson, AZ 85721
\\ email: tgk@math.arizona.edu
}

\maketitle 

\begin{abstract}

We conjecture a relationship between the scaling limit of the fixed-length
ensemble of self-avoiding walks in the upper half plane and radial SLE$_{8/3}$
in this half plane from $0$ to $i$. The relationship is that if we 
take a curve from the fixed-length scaling limit of the SAW, weight it 
by a suitable power of the distance to the endpoint of the curve and 
then apply the conformal map of the half plane that takes the endpoint
to $i$, then we get the same probability measure on curves as radial 
SLE$_{8/3}$. In addition to a non-rigorous derivation of this conjecture, 
we support it with Monte Carlo simulations of the SAW. 
Using the conjectured relationship between the SAW and radial SLE$_{8/3}$, 
our simulations give estimates for 
both the interior and boundary scaling exponents. The values we obtain are 
within a few hundredths of a percent of the conjectured values.

\end{abstract}

\newpage

\section{Introduction}
\label{intro}

Consider the uniform probability measure on all self-avoiding walks 
(SAW) with $N$ steps on a two dimensional lattice with spacing $\delta$, 
e.g., $\delta \integers^2$.  We restrict to SAW's that 
start at the origin and lie in the upper half plane thereafter.
Now take $\delta=N^{-\nu}$ and let $N \ra \infty$. This should 
give a probability measure on simple curves in the half plane that 
start at the origin and end somewhere in the half plane.
We refer to this scaling limit as the fixed-length scaling  limit.
Radial SLE$_{8/3}$ gives a probability measure on simple curves in the 
half plane that start at the origin and end at some prescribed point. 
So it is natural to look for some relation between the 
fixed-length scaling limit of the SAW and radial SLE$_{8/3}$. 

The simplest relation would be the following. Take a curve from the 
fixed-length scaling limit of the half plane SAW 
and apply the conformal map of the half plane 
to itself that fixes the origin and maps the endpoint of the curve to $i$. 
This transformation gives a probability measure on curves from 
the origin to $i$, 
and so one might ask if the resulting measure is radial SLE$_{8/3}$.
Oded Schramm proposed a Monte Carlo study of this possibility 
to the author \cite{schramm_private}.
Simulations of the SAW gave strong evidence that these two 
probability measure are not the same. 
In this paper we will argue that this process of conformally mapping the 
endpoint of the SAW to $i$ does in fact give radial SLE$_{8/3}$ if 
one weights the walk by $R^p$ where $R$ is the distance from the origin 
to the endpoint of the walk taken from the fixed-length scaling limit.
The power $p$ is conjectured to be $(\rho-\gamma)/\nu=-61/48$.

In addition to our heuristic derivation of this conjecture, we also 
provide Monte Carlo simulations of the SAW to support it. 
The Monte Carlo simulations
show excellent agreement with analytic calculations done with 
radial SLE$_{8/3}$. We can also estimate the interior and exterior
scaling exponents from the simulations and find values that 
are within a few hundredths of a percent of the conjectured values. 

One can also consider an ensemble of SAW's in the half plane that 
start at the origin and have any length. The walks are then weighted
by $\mu^{-N}$ where $N$ is the length of the walk and $\mu$ is the 
lattice constant for the SAW on the lattice being used. 
The total weight of these walks is infinite, but 
SLE partition functions predict the relative weights of the walks 
that end at different points in the upper half plane. Our simulations of 
the SAW also provide a partial test of this prediction. More 
precisely, we are able to test the angular dependence of the prediction.
Once again excellent agreement is found. 

In section \ref{background}
we quickly review several different 
scaling limits one can take for the SAW and their conjectured relationships 
with chordal and radial SLE$_{8/3}$. We also review several critical 
exponents for the SAW and review SLE partition function 
predictions for the SAW. 
Section  \ref{conjecture_sect}
is devoted to a non-rigorous derivation of our conjectured relationship
between the fixed-length scaling limit of the SAW and radial SLE$_{8/3}$. 
We give two explicit conjectures about this relationship. One of 
them provides a partial test of SLE partition function predictions 
for the SAW. In section \ref{exact_radial}
we define four random variables that we 
use to test our conjectured relationship. The distributions of these random
variables may be explicitly computed for radial SLE$_{8/3}$ using a 
theorem of Lawler, Schramm and Werner \cite{lsw_restriction}, and we 
give the results of those calculations.
Section \ref{simulations}
gives the results of our Monte Carlo simulations of the 
SAW and the comparison with the radial SLE$_{8/3}$ computations. 
Finally, in section \ref{conclusions} we give our conclusions.

\section{Background}
\label{background}

In this section we review a variety of results, mostly non-rigorous,
on the two-dimensional self-avoiding walk (SAW) that we will use.
Some of the results we review involve the Schramm Loewner Evolution (SLE)
introduced in \cite{schramm}. 
Our conjectures and simulations will only involve ensembles defined in 
the half plane, but in this section we will consider other domains
as well. We refer the reader to one of the 
reviews \cite{bb_review, cardy_review, kn_review, werner_review} 
or the book \cite{lawler_book} for background on SLE. 
A good general reference for the SAW is \cite{ms_book}.

A SAW is a nearest neighbor walk 
on a lattice that does not visit any site more than once.
We denote a SAW by $\omega$ and let $|\omega|$ denote the number
of steps in the walk. So for $i=0,1,\cdots,|\omega|$, 
$\omega(i)$ is a site in the lattice.
For $i=1,\cdots,|\omega|$ we have $||\omega(i)-\omega(i-1)||=1$, and 
$\omega(i) \ne \omega(j)$ for $i \ne j$. 
Our simulations will only be for
SAW's on the square lattice, but our conjecture should hold for other
two-dimensional lattices as well, e.g., 
the triangular and hexagonal lattices. 
There are a variety of scaling limits that one can consider. 
We will first consider a scaling limit which uses SAW's with the same 
number of steps. This is the scaling limit which is most relevant 
from the physics viewpoint and has been studied extensively. 
Unfortunately it is not directly described 
by SLE$_{8/3}$. Next we will consider two scaling limits 
which are conjectured to be directly related to SLE$_{8/3}$. 

If the lattice has unit spacing, the mean square distance 
traveled by an $N$ step SAW grows with $N$ as 
\beann
\EN \, [|\omega(N)|^2] \sim N^{2 \nu}
\eeann
where $|\omega(N)|$ denotes the length of the vector $\omega(N)$ and 
$\EN$ denotes expectation with respect to the uniform probability 
measure on the set of SAW's with $N$ steps that start at the origin. 
The conjectured value of $\nu$ is $3/4$ \cite{flory}.
Now suppose we replace the lattice with unit spacing with a lattice
with spacing $\delta=N^{-\nu}$, and let $N \ra \infty$.  
This scaling limit is expected to give a probability measure on 
simple curves in the plane which start at the origin and end at 
a distance from the origin that is typically of order $1$.

We can also consider this scaling limit for SAW's restricted to the 
upper half plane, $\half$.
We take SAW's with $N$ steps that start at the origin and then stay 
in $\half$ with the uniform probability measure. 
We then take a scaling limit just as we did above by taking the lattice spacing
to be $N^{-\nu}$.
(It is expected that the exponent $\nu$ is the same.) 
This is the SAW ensemble that is the focus of this paper, and we will 
refer to this scaling limit simply as the fixed-length scaling limit. 
We will denote the probability 
measure on curves in $\half$ that start at the origin that comes from 
this scaling limit by $\probfl$ and the corresponding expectation by $\Efl$.

The next scaling limits we consider involve the 
lattice constant $\mu$ defined as follows. 
Let $c_N$ be the number of self-avoiding walks in the full 
plane starting at the origin with $N$ steps. 
It is known that the limit $\lim_{N \ra \infty} c_N^{1/N}$
exists. We denote it by $\mu$. Its value depends on the lattice. 
Nienhuis \cite{nienhuis} conjectured that for the hexagonal lattice
$\mu=\sqrt{2+\sqrt{2}}$. This conjecture was recently proven by 
Duminil-Copin and Smirnov \cite{dc_smirnov}. 
For the square and triangular lattices we have only numerical estimates of 
the value of $\mu$. 

The next two scaling limits we consider are conjectured to 
be directly related to SLE$_{8/3}$. 
We refer to them as the ``chordal scaling limit'' and the 
``radial scaling limit'' since they will correspond to chordal and 
radial SLE$_{8/3}$. 
Consider a simply connected domain $D$. Introduce a lattice with 
spacing $\delta$. Fix a point $z$ on the boundary of $D$ and a 
point $v$ in the interior of $D$. We let $[z]$ and $[v]$ denote the 
lattice sites closest to $z$ and $v$. Then we consider 
all SAW's from $[z]$ to $[v]$ that stay in 
$D$. The number of steps in the SAW's is not constrained. 
We weight each walk $\omega$ by $\mu^{-|\omega|}$ where $|\omega|$ is the 
number of steps in $\omega$. 
The total weight is then 
\bea
Z(D,z,v,\delta) =\sum_{\omega: z \ra v, \omega \subset D} \mu^{-|\omega|}
\label{norm}
\eea
We define a probability measure by weighting each walk $\omega$ by
$\mu^{-|\omega|}/Z(D,z,v,\delta)$.
The scaling limit $\delta \ra 0$ 
is believed to exist and equal radial SLE$_{8/3}$ in $D$ from $z$ to 
$v$ \cite{lsw_saw}. We are primarily interested in the case that 
$D=\half$ and $z=0$. We will use $\prob^{radial}_v$ to denote the 
radial SLE$_{8/3}$ probability measure on curves in $\half$ from $0$ to $v$. 

The chordal scaling limit is similar. The only difference 
is that both $z$ and $v$ are boundary points for $D$. 
It is conjectured that this scaling limit 
is chordal SLE$_{8/3}$ \cite{lsw_saw}. 
If $D=\half$ (or another unbounded domain) 
then one possible choice for $z$ or $v$ is the boundary 
point $\infty$, but we cannot 
use the above construction in this case. Instead we can do the 
following. Use the uniform probability measure on 
SAW's in $\half$ starting at $0$ with $N$ steps.
We first take the limit $N \ra \infty$ to 
get a probability measure on infinite length SAW's. 
Then we let the lattice spacing go to zero. The result is 
believed to be a probability measure on simple curves starting 
at the origin and going to infinity that agrees with chordal SLE$_{8/3}$ 
from the origin to infinity in $\half$. This conjecture was tested
by Monte Carlo simulations of the SAW, and excellent agreement was 
found \cite{Kennedya,Kennedyb}.
We can take this same double limit for SAW's that start 
at the origin and live in the full plane, but this scaling limit will
not play a role in this paper.  

The existence of the limit of the uniform probability measure 
on half-plane SAW's with $N$ steps as $N \ra \infty$ 
has been proved. 
Madras and Slade \cite{ms_book}, using results of Kesten
\cite{kestena,kestenb},
proved that the uniform probability measure on $N$-step bridges 
has a weak limit as $N \ra \infty$. 
Lawler, Schramm and Werner used their methods to prove 
that the uniform probability measure on $N$-step SAW's in the half-plane
has a weak limit as $N \ra \infty$ which is the same as this 
infinite bridge measure \cite{lsw_saw}.  
However, the existence of the scaling limit has not been proved. 

There is another way to construct the probability measure on infinite 
length SAW's in $\half$ that is closer in spirit 
to the definition for bounded domains. Consider all finite length 
SAW's in $\half$ that start at $0$. If we weight a walk 
$\omega$ by $\mu^{-|\omega|}$, then the total weight of the walks will 
be infinite. So we take $x > \mu$ and weight $\omega$ by $x^{-|\omega|}$. 
Then the total weight is finite, and so we can normalize to get a probability
measure. The limit as $x \ra \mu^{+}$ has been proved to exist and 
give the same measure on infinite walks in $\half$ as the probability 
measure one obtains as the limit of the uniform measure on $N$-step SAW's 
as $N \ra \infty$ \cite{saw_bridges}.

Finally we consider how the normalization factor \reff{norm}
depends on $z$ and $v$. 
It is conjectured that there are scaling exponents $b$ and $\bp$ 
and a function $H(D,z,v)$ such that as the lattice spacing goes to zero, 
\bea
Z(D,z,v,\delta) \sim \delta^{b+\bp} H(D,z,v)
\label{sle_pf}
\eea
and $H(D,z,v)$ satisfies the following form of conformal covariance. 
If $\Phi$ is a conformal map of $D$ onto $D^\prime$, $z^\prime=\Phi(z)$ and 
$v^\prime=\Phi(v)$, then 
\bea
H(D,z,v,\delta) = |\Phi^\prime(z)^b \Phi^\prime(v)^{\bp}| \, 
H(D^\prime,z^\prime,v^\prime)
\label{conf_covar} 
\eea
See \cite{lawler_part_func,lawler_utah,lsw_saw}.
We will refer to the function $H(D,z,v)$ as an SLE partition function. 
Note that in \cite{lsw_saw}, the boundary scaling exponent $b$ is denoted
by $a$, and the interior scaling exponent $\bp$ is denoted by $b$. 

For a general domain $D$, equations \reff{sle_pf} and \reff{conf_covar} 
are not the full story. There are expected to be lattice 
effects associated with the boundary point $z$ which 
persist in the scaling limit.
In this paper we only use this equation for the domain $\half$. In 
this case there are no such lattice effects to worry about.

Eq. \reff{conf_covar} determines $H(D,z,v,\delta)$ up to an overall constant.
In particular we can compute $H(\half,0,v)$. It will be convenient to 
represent the endpoint in polar coordinates $r e^{i \theta}$. 
The conformal automorphism of $\half$ that fixes $0$ and maps 
$r e^{i \theta}$ to $i$ is given by 
\bea
\Phi(z)= {z \sin \theta \over r- z \cos \theta}
\eea
\ifodd \verbose
\bea
\Phi^\prime(z)= {r \sin \theta \over (r- z \cos \theta)^2}
\eea
From which we get 
\beann
|\Phi^\prime(0)| &=& {\sin \theta \over r}, 
\nonumber \\
|\Phi^\prime(r e^{i \theta})| &=& {1 \over r \sin \theta } 
\eeann
\else
\fi
With the convention that $H(\half,0,i)=1$ we find 
\bea
H(\half,0,r e^{i \theta}) = r^{-b-\bp} (\sin \theta)^{b - \bp}
\label{explicit_sle_pf}
\eea

There are many critical exponents associated with the SAW. 
In addition to the exponent $\nu$ we will need two others:
$\gamma$ and $\rho$. 
The definition of $\mu$ means that the geometric growth of the 
number of SAW's is given by $c_N \sim \mu^N$. 
It is expected that there is a power law
correction to this growth that is usually expressed in the form 
\beann
c_N \sim N^{\gamma-1} \mu^N 
\eeann
The conjectured value of $\gamma$ is $43/32$ \cite{nienhuis}.

Let $b_N$ be the number 
of SAW's that start at the origin and then stay in the upper half plane. 
Restricting the SAW to stay in the half plane does not change the 
geometric rate of growth, but it does change the power law correction. 
It can be written in the form 
\bea
b_N \sim N^{\gamma-1-\rho} \mu^N 
\label{rhodef}
\eea
Note that $b_N/c_N$ gives the probability that a full-plane 
SAW with $N$ steps stays in 
the upper half plane. It goes as $N^{-\rho}$, so the exponent $\rho$ 
characterizes this probability. It is conjectured that $\rho=25/64$
\cite{lsw_saw}.

\section{The conjecture}
\label{conjecture_sect}

To state our conjecture precisely we introduce some notation. 
We use $\probNdelta$ to denote the uniform probability measure on $N$-step
SAW's in $\half$ on the lattice with spacing $\delta$ that start at the 
origin. We use $\ENdelta$ to denote the corresponding expectation.
We assume that if we take $\delta=N^{-\nu}$ and let $N \ra \infty$, then 
the scaling limit of $\probNdelta$ exists. 
We denote it by $\probfl$ and refer to it as the  fixed-length 
scaling limit of the SAW. The corresponding expectation is denoted $\Efl$.
We use $\probsle$ to denote the probability measure for radial SLE 
in $\half$ from $0$ to $i$.

\bigskip

\no {\bf \large Conjecture 1:} 
The fixed length scaling limit of the SAW
and radial SLE in the half plane from $0$ to $i$ are related by 
\bea
{\Efl[R(\gamma)^{(\rho-\gamma)/\nu} \, 1(\phi_\gamma(\gamma) \in E)]
\over \Efl[R(\gamma)^{(\rho-\gamma)/\nu}]}
= \probsle (E).
\label{conjectureone}
\eea
$E$ is an event for simple curves in $\half$ that go from $0$ to $i$. 
In the left side, $\gamma$ is the random curve from the fixed length 
scaling limit measure, 
$R(\gamma)$ is the distance from the origin to the endpoint of 
$\gamma$, and $\phi_\gamma$ is the Moibius transformation of $\half$ that 
fixes the origin and takes the endpoint of $\gamma$ to $i$. 
So we can generate radial SLE$_{8/3}$ by generating a 
curve $\gamma$ from the fixed-length scaling limit, weighting it by 
$R(\gamma)^{(\rho-\gamma)/\nu}$, and applying the conformal map that takes 
the endpoint to $i$. 
With the conjectured values of the exponents, 
the power $(\rho-\gamma)/\nu$ is $-61/48$.

\bigskip

We give a non-rigorous derivation of this conjecture 
in this section.
We consider SAW's on a lattice with spacing $\delta$. 
Fix $0<r_1<r_2$. (We think of them as being of order 1. They will 
not diverge or go to zero in the scaling limit.) 
We consider all SAW's in $\half$ which start at the origin and 
end somewhere in the region 
\beann
A=\{ z \in \half : r_1 \le |z| \le r_2\}.
\eeann
There is no constraint on the number of steps in the SAW. 
We weight a SAW $\omega$ by $\mu^{-|\omega|}$. 
(As before $|\omega|$ is the number of steps in $\omega$.)
The total weight of the walks that end in the annular region is finite. 
(This is the reason for introducing the cutoffs $r_1$ and $r_2$.)
So we can normalize to obtain a probability measure. 
To be more precise, we define
\bea
Z(A,\delta) = \sum_{\omega} \mu^{-|\omega|}
1(R(\omega) \in [r_1,r_2]) \,
\eea
where $R(\omega)$ is the distance of the endpoint of $\omega$ to the origin.
The sum over $\omega$ is over all SAW's in $\half$ that start at $0$.
The constraint that $\omega$ ends in $A$ is incorporated 
in the indicator function.
In the notation of the previous section, $Z(\half,0,v,\delta)$ would 
denote the weight of the walks that end at $v$. So it would have been 
more consistent to denote the above by 
$Z(\half,0,A,\delta)$. Since all SAW's in this section start at $0$ 
and stay in $\half$, we have shortened this to just $Z(A,\delta)$. 
We then assign probability $\mu^{-|\omega|}/Z(A,\delta)$ to each 
walk that ends in the annular region.

The above probability measure is on SAW's that end in $A$. We now use it 
to define a probability measure on curves in $\half$ that go from $0$
to $i$. 
For a SAW $\omega$, let $\phi_\omega(z)$ be the conformal automorphism 
of $\half$ that fixes $0$ and takes the endpoint of $\omega$ to $i$. 
(Of course, this map only depends on the endpoint of $\omega$, not the 
entire SAW.) The image $\phi_\omega(\omega)$ will be a curve in $\half$ 
from $0$ to $i$. The probability measure of the previous paragraph
gives a probability measure on such curves. 
Let $E$ be an event for such curves. 
The probability of $E$ 
is defined to be $N(E,A,\delta)/Z(A,\delta)$, where 
\bea
N(E,A,\delta) = \sum_{\omega} \mu^{-|\omega|}
1(R(\omega) \in [r_1,r_2]) \, 1(\phi_\omega(\omega) \in E)
\label{initial_numer}
\eea

We decompose the sum over walks by length. 
\beann
N(E,A,\delta) &=& \sum_N \mu^{-N} 
\sum_{\omega:|\omega|=N} \, 1(R(\omega) \in [r_1,r_2]) 
\, 1(\phi_\omega(\omega) \in E)
\nonumber \\
&=& \sum_N \mu^{-N} b_N \, \, \ENdelta [ 1(R(\omega) \in [r_1,r_2]) 
\, 1(\phi_\omega(\omega) \in E)]
\nonumber \\
\eeann
Using \reff{rhodef} we replace $\mu^{-N} b_N$ by $N^{\gamma-1-\rho}$.
\beann
N(E,A,\delta) 
&\approx& \sum_N N^{\gamma-1-\rho} 
\, \, \ENdelta [ 1(R(\omega) \in [r_1,r_2])
\, 1(\phi_\omega(\omega) \in E) ]
\nonumber \\
\eeann
When the lattice spacing $\delta$ is small, the constraint that 
the SAW ends at a distance from the origin that lies in $[r_1,r_2]$
implies that $N$ must be large. So we can approximate $\ENdelta$ using the
fixed-length scaling limit. If we rescale $\omega$ by a factor
of $\delta^{-1} N^{-\nu}$, then in the limit its
distribution will converge to that of $\gamma$ drawn from $\probfl$,
the probability measure of the fixed-length scaling limit. 
We rewrite the condition that $R(\omega) \in [r_1,r_2]$ as 
$r_1 N^{-\nu} \delta^{-1} \le R(\omega) N^{-\nu} \delta^{-1}
\le r_2 N^{-\nu} \delta^{-1}$, 
so this becomes the condition
$r_1 N^{-\nu} \delta^{-1} \le R(\gamma) \le r_2 N^{-\nu} \delta^{-1}$
where $R(\gamma)$ is the distance of the endpoint of $\gamma$ from the 
origin. 
Note that the condition $\phi_\omega(\omega) \in E$ just becomes 
$\phi_\gamma(\gamma) \in E$. 

We now have 
\bea
N(E,A,\delta) 
&\approx& \sum_N N^{\gamma-1-\rho} 
\, \, \Efl [ 1(r_1 N^{-\nu} \delta^{-1} \le R(\gamma) \le
r_2 N^{-\nu} \delta^{-1})
\, \, 1(\phi_\gamma(\gamma) \in E) ]
\eea
We move the sum on $N$ inside the expectation and then consider
\beann
&& \sum_N N^{\gamma-1-\rho} 
\, 1(r_1 N^{-\nu} \delta^{-1} \le R(\gamma) \le r_2 N^{-\nu} \delta^{-1})
\nonumber \\
&=& 
\sum_N N^{\gamma-1-\rho} 
\, 1\left( \left({1 \over \delta} {r_1 \over R(\gamma)}\right)^{1/\nu}  
\le N \le
\left( {1 \over \delta} {r_2 \over R(\gamma)}\right)^{1/\nu} \right)
\eeann

Since $\delta \ra 0$, the values of $N$ are large, and so we can replace
the above by 
\beann
&&\int_0^\infty \, x^{\gamma-1-\rho} 
\, 1\left( \left({1 \over \delta} {r_1 \over R(\gamma)}\right)^{1/\nu}  
\le x \le
\left( {1 \over \delta} {r_2 \over R(\gamma)}\right)^{1/\nu} \right) \, dx
\nonumber \\
&=& \delta^{(\rho-\gamma)/\nu} \, R(\gamma)^{(\rho-\gamma)/\nu} \,
\int_{r_1^{1/\nu}}^{r_2^{1/\nu}} \, x^{\gamma-1-\rho} \, dx
\nonumber \\
&=& 
c \, \delta^{(\rho-\gamma)/\nu} \, R(\gamma)^{(\rho-\gamma)/\nu} 
\label{nsum}
\eeann
The factor of $c \delta^{(\rho-\gamma)/\nu}$ will cancel with the corresponding 
factor in $Z(A,\delta)$. 
Note that the approximation in the above corresponds to the rigorous 
statement that 
\bea
\lim_{\delta \ra 0} \delta^{(\gamma-\rho)/\nu} 
\, \sum_N N^{\gamma-1-\rho} 
\, 1\left(\left({1 \over \delta} {r_1 \over R(\gamma)}\right)^{1/\nu}  
\le N \le
\left({1 \over \delta} {r_2 \over R(\gamma)}\right)^{1/\nu} \right)
= c \, R(\gamma)^{(\rho-\gamma)/\nu} 
\eea
So we find
\bea
\lim_{\delta \ra 0} {N(E,A,\delta) \over Z(A,\delta)}
= {\Efl[R(\gamma)^{(\rho-\gamma)/\nu} \, 1(\phi_\gamma(\gamma) \in E)]
\over \Efl[R(\gamma)^{(\rho-\gamma)/\nu}]}
\eea

We now return to \reff{initial_numer} and decompose the sum 
according to the endpoint of the walk. 
\bea
N(E,A,\delta) = \sum_{z \in \delta \integers^2 \cap A} \, 
\sum_{\omega : 0 \ra z} \mu^{-|\omega|} \, 1(\phi_\omega(\omega) \in E)
\eea
where the notation $\omega : 0 \ra z$ means that the sum over $\omega$ is 
over walks between $0$ and $z$. Recall that
\bea
Z(\half,0,z,\delta) = \sum_{\omega: 0 \ra z} \mu^{-|\omega|} 
\eea
Let $\probzdelta$ denotes the corresponding probability measure on SAW's from
$0$ to $z$ which gives $\omega$ the weight $\mu^{-|\omega|}/Z(\half,0,z,\delta)$. 
Then we can rewrite $N(E,A,\delta)$ as
\bea
N(E,A,\delta) = \sum_{z \in \delta \integers^2 \cap A}
Z(\half,0,z,\delta) \, \probzdelta (\phi_\omega(\omega) \in E)
\eea
As $\delta \ra 0$, $\probzdelta (\phi_\omega(\omega) \in E)$ should 
converge to $\probsle (E)$, where $\probsle$ denotes the 
probability measure for radial SLE$_{8/3}$ in $\half$ from $0$ to $i$. 
So
\bea
\lim_{\delta \ra 0} {N(E,A,\delta) \over Z(A,\delta)}
= \probsle (E).
\eea
Thus we have derived the conjecture \reff{conjectureone}.

Recall that the scaling limit of 
$Z(\half,0,z,\delta)$ is conjectured to be given by \reff{sle_pf} with 
$H(\half,0,z)$ given by \reff{explicit_sle_pf}.
This formula did not enter the derivation
of conjecture \reff{conjectureone}. We will derive a second conjecture
that does involve, at least partially, this SLE partition function.
For an angle $\theta \in [0,\pi]$ we define
\bea
F(\theta,A,\delta) = \sum_{\omega} \mu^{-|\omega|}
1(r_1 \le R(\omega) \le r_2) \, 1(arg(\omega) \le \theta)
\label{initial_f}
\eea
where $arg(\omega)$ denotes the polar angle of the endpoint of the SAW.
The same argument as before shows 
\beann 
\lim_{\delta \ra 0} \delta^{(\gamma-\rho)/\nu} \, F(\theta,A,\delta) 
&=& c \, \Efl[R(\gamma)^{(\rho-\gamma)/\nu} \, 1(\arg(\gamma) \le \theta)]
\eeann

If we decompose the sum in \reff{initial_f} 
according to the endpoint of the walk we have 
\bea
F(\theta,A,\delta) &=& 
\sum_{z \in \delta \integers^2 \cap A}
1(arg(z) \le \theta) \, \sum_{\omega:0 \ra z} \mu^{-|\omega|}
\nonumber \\
&=& \sum_{z \in \delta \integers^2 \cap A}
1(arg(z) \le \theta) \, Z(\half,0,z,\delta) 
\eea
As $\delta \ra 0$, $ \delta^{-b-\bp} Z(\half,0,z,\delta)$ 
should converge to the SLE partition function $H(\half,0,z)$. 
We assume that 
\bea
b+\bp= {\rho-\gamma \over \nu} +2 
\eea
This relation was conjectured in \cite{lsw_saw}. It is satisfied by 
the conjectured values $b=5/8$, $\bp=5/48$, $\gamma=43/32$ and $\rho=25/64$.
So we can rewrite the above as 
\bea
\delta^{(\gamma-\rho)/\nu} F(\theta,A,\delta) 
= \delta^2 \sum_{z \in \delta \integers^2 \cap A}
1(arg(z) \le \theta) \, \delta^{-b-\bp} Z(\half,0,z,\delta) 
\eea
As $\delta \ra 0$, the sum over $z$ together with the factor of 
$\delta^2$ becomes an integral over $z$.
Then using \reff{explicit_sle_pf} we obtain our second conjecture. 

\bigskip

\no {\bf \large Conjecture 2:} 
\bea
{\Efl \left[R(\gamma)^{(\rho-\gamma)/\nu} \, 1(\arg(\gamma) \le \theta) \right]
\over \Efl \left[R(\gamma)^{(\rho-\gamma)/\nu} \right]}
= { \int_0^\theta \, \sin(\alpha)^{b-\bp} d \alpha
\over \int_0^\pi \, \sin(\alpha)^{b-\bp} d \alpha}
\label{conjecturetwo}
\eea

\section{Exact calculations for SLE}
\label{exact_radial}

In this section we compute the distribution of four random variables
defined in terms of radial SLE$_{8/3}$ in $\half$ from $0$ to $i$. 
We will use these and simulations of 
the SAW to test our conjecture. We will also use 
these distributions and the simulations to estimate the values
of the scaling exponents $b$ and $\bp$ and compare them to the 
conjectured values. 

We consider radial SLE$_{8/3}$ in the half plane from $0$ to $i$.
We denote the probability measure by $\probsle$ and 
the SLE curve by $\gamma$.
Suppose $A$ is a closed set not containing $0$ or $i$, such that 
$\half \setminus A$ is simply connected. We want to compute the 
probability that the SLE curve does not enter $A$. 
Let $\phi_A$ be the 
conformal map from $\half \setminus A$ onto $\half$ with
$\phi_A(0) = 0$ and $\phi_A(i)=i$. Then 
\bea
\probsle(\gamma \cap A = \emptyset) = 
|\phi_A^\prime(0)|^b \, |\phi_A^\prime(i)|^{\bp}, \quad
b=5/8, \, \bp =5/48
\label{lsw_formula}
\eea 
This formula is stated at the end of \cite{lsw_restriction}. 
(They state a formula for radial SLE$_{8/3}$ in the unit disc which 
immediately gives the above formula.) 

\ifodd \verbose
The precise formula in \cite{lsw_restriction} is as follows.  
Consider radial SLE$_{8/3}$ in the unit disc $\disc$ from 1 to 0. 
If $\gamma$ is a radial path in $\disc$ from 1 to 0, and $A$
is a compact set not containing 1, such that $\disc \setminus A$
is simply connected and contains 0, and $\psi_A$ is 
a conformal map from $\disc \setminus A$ onto $\disc$ with
$\psi_A(0) = 0$, then
\bea
\probsle(\gamma[0,\infty) \cap A = \emptyset) = 
|\psi_A^\prime(1)|^b \, |\psi_A^\prime(0)|^{\bp} 
\label{lsw_restrictionb}
\eea 
To derive \reff{lsw_formula} from this let $\phi$ be the conformal 
map of $\half$ to $\disc$ with $\phi(0)=i$ and $\phi(0)=1$. 
The probability that SLE in $\half$ misses $A$ equals 
the probability that SLE in $\disc$ misses $\phi(A)$.
Let $\psi=\phi \compose \phi_A \compose \phi^{-1}$. 
Then $\psi$ maps $\disc \setminus \phi(A)$ to $\disc$, and fixes 
$0$ and $1$. So \reff{lsw_restrictionb} gives 
\bea
\prob^{radial}_\disc(\gamma[0,\infty) \cap \phi(A) = \emptyset) = 
|\psi^\prime(1)|^{b} \, |\psi^\prime(0)|^{\bp} 
\eea 
A chain rule calculation shows
$\psi^\prime(0)=\phi_A^\prime(i)$ and 
$\psi^\prime(1)=\phi_A^\prime(0)$. So \reff{lsw_formula} follows.
\else
\fi

The first random variable we consider is the rightmost excursion of 
$\gamma$. So
\bea
X = \max_t Re (\gamma(t))
\eea
Since $\gamma$ starts at the origin, $X \ge 0$. 
We want to compute the probability that $X < x$, i.e., that the 
walk does not entered the region given by the quarter plane 
$A_x = \{z \in \half : Re(z) \ge x\}$.   
\ifodd \verbose
Let $f_x(z)=x^2 -(z-x)^2$. Then $f$ maps $\half \setminus A_x$ 
to $\half$ and $f(0)=0$. We have $f(i)=1+2ix$. 
Let 
\bea 
\psi_x(z)={2zx \over 4x^2 +1 -z}
\eea
This is the automorphism of $\half$ that fixes $0$ and sends
$1+2ix$ to $i$. 
The conformal map $\phi_{A_x}$ is then given by 
$\phi_x(z)=\psi_x(f_x(z))$.
\else
\fi
We denote the conformal map $\phi_{A_x}$ by just $\phi_x$. It is given by 
\bea
\phi_x(z)= {2x (2xz-z^2) \over z^2 -2xz+4x^2+1}
\eea
\ifodd \verbose
\bea
\phi^\prime_x(z)=
{4x (4x^2+1)(x-z) \over (z^2 -2xz+4x^2+1)^2}
\eea
And so
\beann
\phi^\prime_x(0) &=& {4x^2 \over 4x^2+1}
\nonumber \\
\phi^\prime_x(i) &=& { (4x^2+1)(x-i) \over  x (2x-i)^2}
\nonumber \\
|\phi^\prime_x(i)| &=& { \sqrt{x^2+1} \over  x }
\eeann
\else
\fi
After some computation, \reff{lsw_formula} yields
\bea
\probsle (X \le x) = 
\left[{4x^2 \over 4x^2+1}\right]^{b} 
\left[{\sqrt{x^2+1} \over  x }\right]^{\bp}
\eea

The second random variable we consider is the highest excursion of 
$\gamma$. So
\bea
Y = \max_t Im (\gamma(t))
\eea
The event $Y < y$ says that the walk does not enter the 
half plane $A_y=\{z : Im(z) \ge y\}$. 
\ifodd \verbose
So we need the conformal map $\phi_y(z)$ from the strip 
$\{z : 0 < Im(z) < y\}$ to $\half$ with $\phi_y(0)=0$ and 
$\phi_y(i)=i$. Let 
\bea
f(z)=\exp({\pi z \over y}) -1 
\eea
This maps the strip to $\half$ and fixes $0$ but not $i$. 
The image of $i$ is $\exp({i \pi \over y}) -1$. 
Let 
\bea
\psi(z)= {\sin({\pi \over y}) z \over 
(1-\cos({\pi \over y})) (z+2) }
\eea
This is the conformal automorphism of $\half$ which fixes $0$ and 
sends $\exp({i \pi \over y}) -1$ to $i$. 
So $\phi_y = \psi \compose f$.  
This is given by 
\bea
\phi_y(z) &=& {\sin({\pi \over y}) [\exp({\pi z \over y})-1] 
\over [1-\cos({\pi \over y})] [\exp({\pi z \over y}) + 1]}
\eea
\else
\fi
The conformal map $\phi_{A_y}$ is given by 
\bea
\phi_y(z) = {\sin({\pi \over y}) \over 1-\cos({\pi \over y})}
\tanh({\pi z \over 2 y})
\eea
\ifodd \verbose
So 
\bea
\phi^\prime_y(z) &=& {\sin({\pi \over y}) \over 1-\cos({\pi \over y})}
\, {\pi \over 2 y} \, \sech^2({\pi z \over 2 y})
\nonumber \\ 
\phi^\prime_y(0) &=& {\sin({\pi \over y}) \over 1-\cos({\pi \over y})}
\, {\pi \over 2 y} 
\nonumber \\ 
\phi^\prime_y(i) &=& {\sin({\pi \over y}) \over 1-\cos({\pi \over y})}
\, {\pi \over 2 y} \, \sech^2({i \pi \over 2 y})
\nonumber 
\eea
\else
\fi
We then find from \reff{lsw_formula} that
\bea
\probsle(Y \le y) = 
\left[{\sin({\pi \over y}) \over 1-\cos({\pi \over y})} \, {\pi \over 2 y}
\right]^{b+\bp} \left[\cos({\pi \over 2 y})\right]^{-2 \bp}
\eea

The third random variable is the maximum distance of $\gamma$ from its starting 
point at the origin. So 
\bea
R = \max_t |\gamma(t)|
\eea
Obviously, $R \ge 1$. 
The event $R < r$ corresponds to $\gamma$ not entering 
$A_r= \{ z \in \half : |z| \ge r \}$. 
\ifodd \verbose
Let $f_r(z)= -(z/r+r/z)$. This maps $\half \setminus A_r$ to $\half$, and
$f_r(0)=\infty, f_r(i)=i(r-1/r)$. 
Let 
\bea
\psi_r(z) ={1-r^2 \over rz} 
\eea
This is the automorphism of $\half$ that sends $\infty$ to $0$ and 
$i(r-1/r)$ to $i$. 
The conformal map is $\phi_r(z)=  \psi(f(z))$ 
\else
\fi
The conformal map is 
\bea
\phi_r(z)=  {(r^2-1) z \over z^2+r^2}
\eea
\ifodd \verbose
So
\bea
\phi^\prime_r(z)=  {(r^2-1)x \over (z^2+r^2)^2}
\eea
and this gives
\bea
\phi_r^\prime(0)= {r^2-1 \over r^2}
\nonumber \\
\phi_r^\prime(i)= {r^2+1 \over r^2-1}
\eea
\else
\fi
So we find 
\bea
\probsle(R \le r) = 
\left[{r^2-1 \over r^2}\right]^{b}
\left[{r^2+1 \over r^2-1}\right]^{\bp}
\eea

Our final random variable is the maximum distance of the walk from the 
endpoint at $i$. 
\bea
S = \max_t |\gamma(t)-i|
\eea
Since $\gamma$ starts at $0$, $S \ge 1$. As we will see,
$\probsle(S=1)>0$. 
The event $S < s$ is that $\gamma$ does not enter
$D_s=\{z : |z-i| \ge s\}$.
Define $l>0$ by $1+l^2=s^2$, so that the circle $|z-i|=s$ intersects the 
real axis at $-l$ and $l$.  
Let 
\beann
f(z)={l+z \over l-z}, \quad g(z)=z^{\pi/\theta}-1,
\eeann
Then $f$ maps $\half \setminus D_s$ to a  wedge 
$0 < arg(z) < \theta$, and $g$ maps this wedge to $\half$.  
The point $l+2i$ is on the circle $|z-i|=s$, and $f(l+2i)=-1+li$. 
So $\theta$ is given by $\tan(\theta)= -l$. 

The map $g \compose f$ fixes $0$. 
Let $x+iy=g(f(i))=\exp(i\pi \alpha/\theta)-1$. 
The automorphism of $\half$ that fixes $0$ and 
takes $x+iy$ to $i$ is 
\bea
\psi(z)={yz \over x^2+y^2-xz}
\eea
So the final conformal map is $\phi_s(z)=\psi(g(f(z)))$.

\ifodd \verbose

Keeping in mind that $f(0)=1, g(f(0))=0$, and
$f(i)=v, g(f(i))=w$, we have 
\bea
\phi_d^\prime(0)= \psi^\prime(0) g^\prime(1) f^\prime(0) \nonumber \\
\phi_d^\prime(i)= \psi^\prime(w) g^\prime(v) f^\prime(i) \nonumber \\
\eea

We have 
\beann
f^\prime(z) = {2 l \over (l-z)^2}, \quad f^\prime(0) = {2 \over l}, \quad
|f^\prime(i)| = {2 l \over l^2 +1}
\eeann
\beann
g^\prime(z)= {\pi \over \theta} z^{\pi/\theta -1}, \quad
g^\prime(1)= {\pi \over \theta}, \quad
|g^\prime(v)| = {\pi \over \theta} 
\eeann

Finally, 
\beann
\psi^\prime(z) = {y (x^2 + y^2) \over (x^2 + y^2 -xz)^2}, \quad
\psi^\prime(0) = {y \over x^2 + y^2}, \quad
|\psi^\prime(w)| = {1 \over y} 
\eeann
From $w=x+iy=\exp(i\pi \alpha/\theta)-1$ we have 
\beann
x^2+y^2 = |w|^2 = 2(1-\cos({\pi \alpha \over \theta})), \quad
y = \sin({\pi \alpha \over \theta})
\eeann
So 
\beann
\psi^\prime(0) = {\sin({\pi \alpha \over \theta}) \over 
2(1-\cos({\pi \alpha \over \theta})) }, \quad
|\psi^\prime(w)| = {1 \over \sin({\pi \alpha \over \theta})}
\eeann
\else
\fi

Note that $|f(i)|=1$, so $f(i)=\exp(i \alpha)$ with $\alpha$ 
given by $\tan \alpha = 2 l /(l^2 -1)$. So 
$x+iy=\exp(i \alpha \pi/\theta) -1$.
Computing all the derivatives we find
\bea
\probsle(S \le s) = \left[{\pi \sin({\pi \alpha \over \theta}) \over 
   l \theta (1-\cos({\pi \alpha \over \theta}))  } \right]^{b}
\left[{2 \pi l  \over \theta \sin({\pi \alpha \over \theta})
   (l^2 +1) } \right]^{\bp}
\eea
where $\theta,\alpha$ and $l$ depend on $s$ through 
$1+l^2=s^2$, $\tan \theta= -l$, and $\tan \alpha = 2 l /(l^2 -1)$.
By taking $s=1$ in the above, we find $\probsle(S=1)=2^{\bp-b}$.

\section{Simulations}
\label{simulations}

The pivot algorithm provides a fast Markov chain Monte Carlo algorithm
for simulating the fixed length ensemble of the
SAW in the full plane or the half plane. 
For an introduction to this algorithm see \cite{ms_book}. We use the 
version of the algorithm found in \cite{kennedy_pivot}, but note 
that a much faster version of the algorithm has been developed by 
Clisby \cite{clisby}.

We simulate the SAW in the half plane with four different numbers of steps: \\
$N=100K, 200K, 500K$ and $1000K$. The iterations of the Markov chain are
highly correlated, so there is no point in sampling the chain at every 
iteration. Instead we sample every $100$ iterations. 
The number of samples generated for each of the four values of $N$ are 
given in table \ref{table_exps}.

\begin{figure}[tbh]
\includegraphics{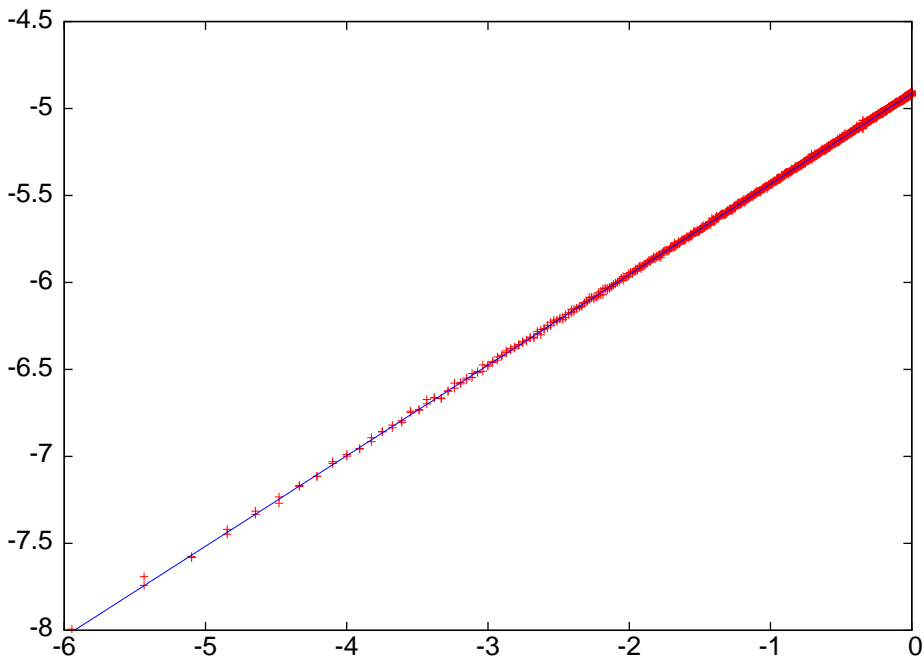}
\caption{\leftskip=25 pt \rightskip= 25 pt 
log-log plot of 
$\EN[R^{(\rho-\gamma)/\nu} \, 1(\theta \le \Theta \le \theta+d \theta)]$
vs. $\sin(\theta+d \theta/2)$. The line is the least 
squares fit to the data.
}
\label{least_sq}
\end{figure}

We first test our second conjecture \reff{conjecturetwo}.
We divide the range of $\Theta$ into $1800$ equal subintervals 
and compute 
$\EN[R^{(\rho-\gamma)/\nu} \, 1(\theta \le \Theta \le \theta+d \theta)]$ 
for each subinterval. We make the approximation
\beann
\int_\theta^{\theta+d\theta} \, \sin(\alpha)^{b-\bp} d \alpha
\approx \left[\sin(\theta+{d\theta \over 2})\right]^{b-\bp} \, d \theta
\eeann
We then do a log-log plot of 
$\EN[R^{(\rho-\gamma)/\nu} \, 1(\theta \le \Theta \le \theta+d \theta)]$
as a function of $\sin(\theta+d \theta/2)$.  
The result for the $N=1000K$ data is 
shown in figure \ref{least_sq}.
The conjecture \reff{conjecturetwo} says that the points should lie on a line. 
The line shown in the figure is a least squares fit to the data. 
It has a slope of \slope. This should be compared with the 
conjectured value of $b-\bp=25/48=0.52083\overline{3}$. 

Next we test our main conjecture \reff{conjectureone}. In our simulations
we generate SAW's from the uniform 
probability measure on SAW's in the half plane with $N$ steps. 
Let $\gamma$ be the SAW scaled by a factor of $N^{-\nu}$.
It is given the weight $R(\gamma)^{(\rho-\gamma)/\rho}$. 
We apply the conformal map that takes the endpoint of the walk to $i$ 
and compute the four random variables $X,Y,R$ and $S$ for the transformed 
walks. Note that we have to normalize by the sum of the weights 
$R(\gamma)^{(\rho-\gamma)/\rho}$, not by the number of samples generated. 
Throughout this section we use $\probNp$ to denote this probability
measure. It depends on the length of the walks that we use in the 
simulation, but for large $N$ it should be a good approximation 
to the ratio in the left side of our first conjecture \reff{conjectureone}.
Our conjecture is that $\probNp$ converges to $\probsle$ as $N \ra \infty$.

\begin{figure}[tbh]
\includegraphics{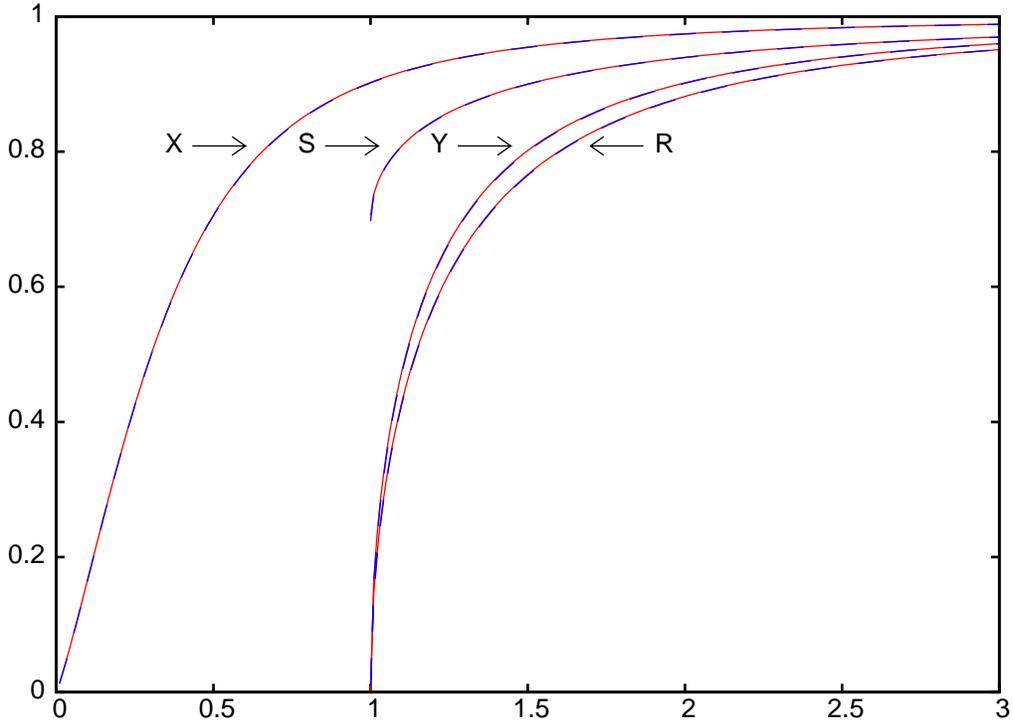}
\caption{\leftskip=25 pt \rightskip= 25 pt 
For each of the random variables $X,Y,R$ and $S$ we plot both the 
cumulative distribution function for radial SLE$_{8/3}$ 
and for the transformed, appropriately weighted SAW. The SLE$_{8/3}$ 
and SAW curves are drawn with dashed lines in different colors. 
They agree so well that the two dashed lines appear to form a single 
curve.
}
\label{rv_all}
\end{figure}

\begin{figure}[tbh]
\includegraphics{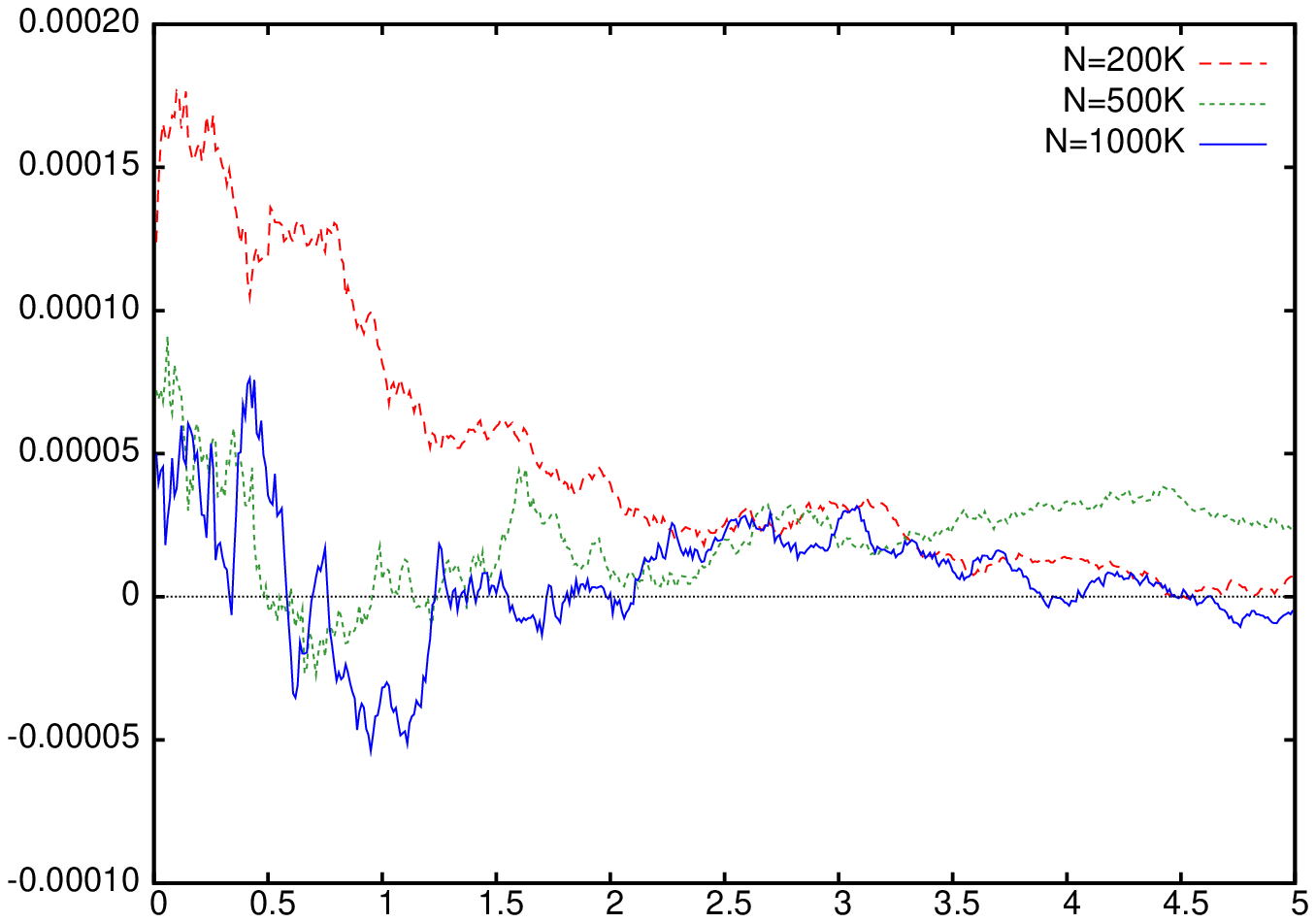}
\caption{\leftskip=25 pt \rightskip= 25 pt 
For the random variable $X$ we plot the difference of the 
cumulative distribution functions for the transformed, 
appropriately weighted SAW and for radial SLE$_{8/3}$. 
}
\label{rv1}
\end{figure}

\begin{figure}[tbh]
\includegraphics{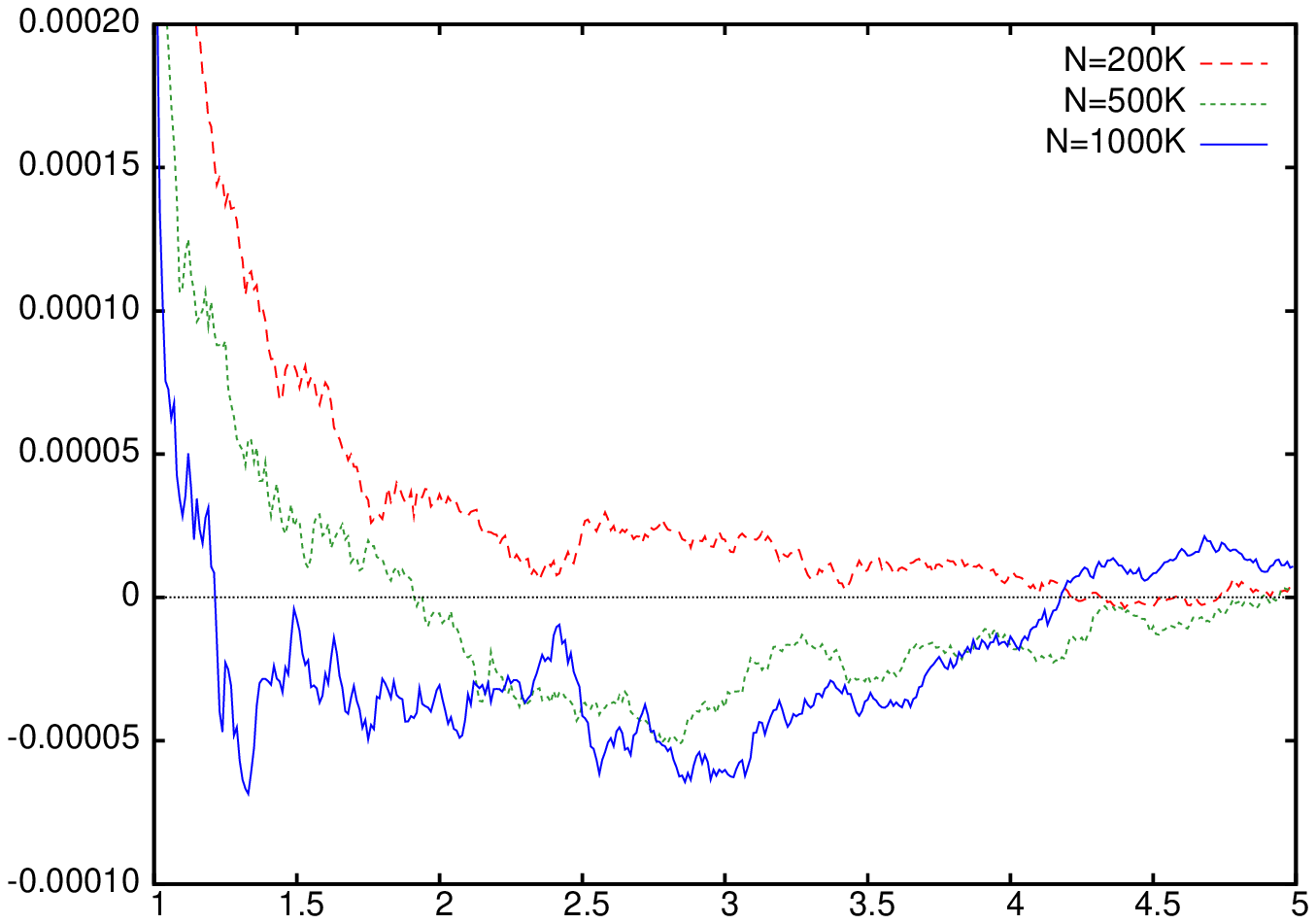}
\caption{\leftskip=25 pt \rightskip= 25 pt 
For the random variable $Y$ we plot the difference of the 
cumulative distribution functions for the transformed, 
appropriately weighted SAW and for radial SLE$_{8/3}$. 
}
\label{rv2}
\end{figure}

\begin{figure}[tbh]
\includegraphics{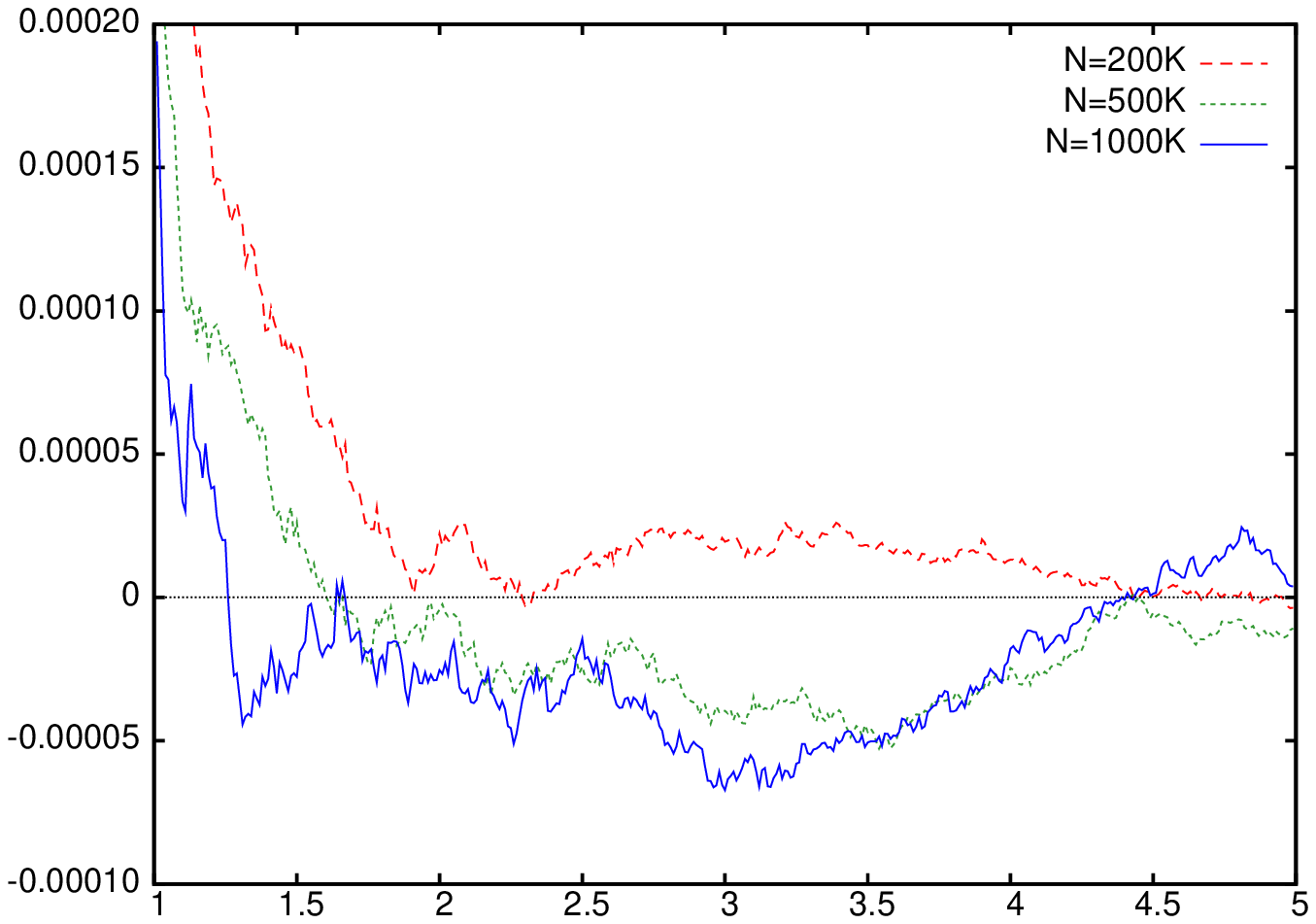}
\caption{\leftskip=25 pt \rightskip= 25 pt 
For the random variable $R$ we plot the difference of the 
cumulative distribution functions for the transformed, 
appropriately weighted SAW and for radial SLE$_{8/3}$. 
}
\label{rv3}
\end{figure}

\begin{figure}[tbh]
\includegraphics{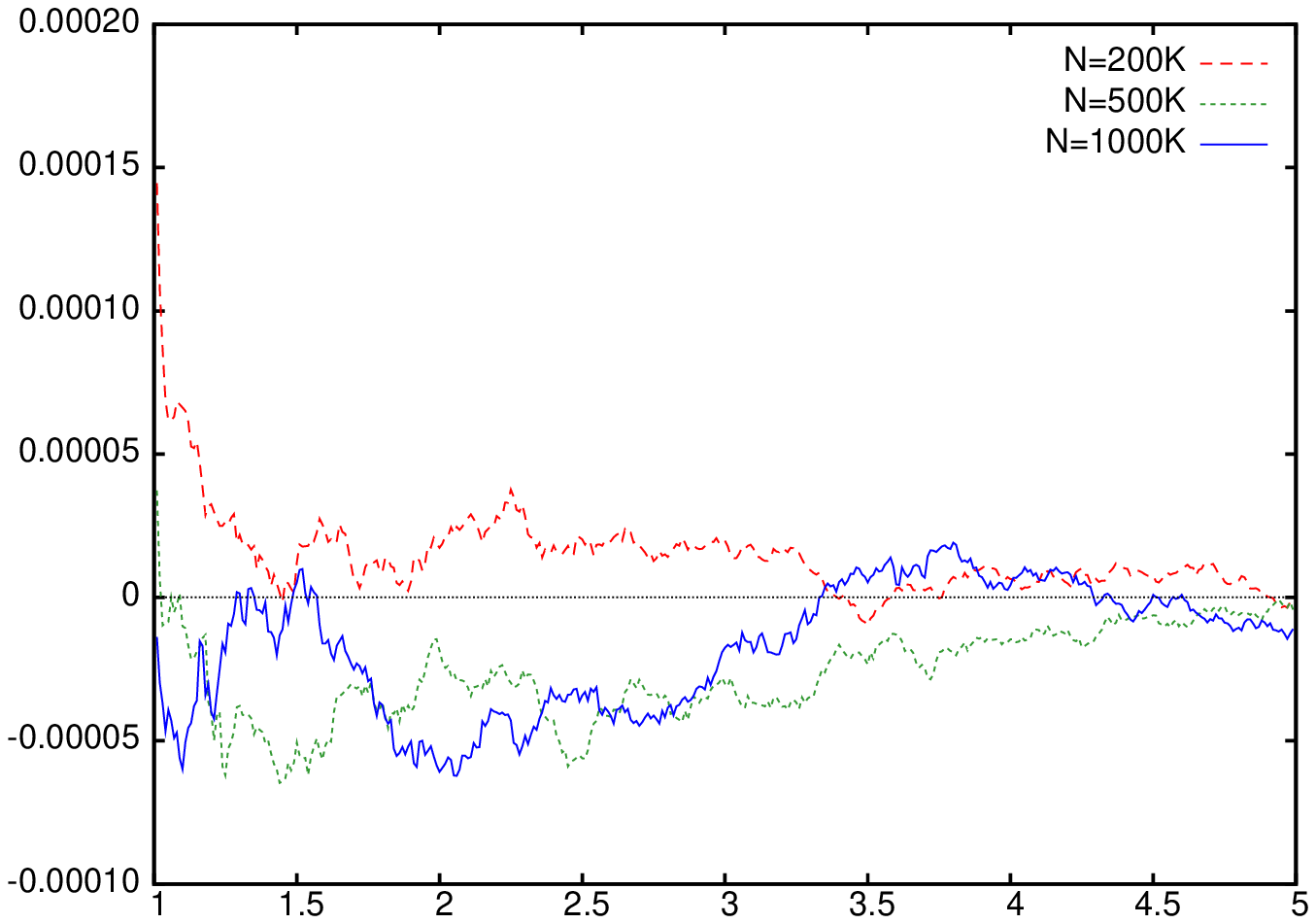}
\caption{\leftskip=25 pt \rightskip= 25 pt 
For the random variable $S$ we plot the difference of the 
cumulative distribution functions for the transformed, 
appropriately weighted SAW and for radial SLE$_{8/3}$. 
}
\label{rv4}
\end{figure}

The cumulative distribution functions from our simulations are shown in 
figure \ref{rv_all} along with the exact  distributions of these 
random variables for radial SLE$_{8/3}$ 
that we computed in section \ref{exact_radial}.
There are eight curves in this figure, but the 
differences between the exact radial SLE$_{8/3}$ results and the simulation
are too small to be seen in the figure so it appears there 
are only four curves. The differences between the analytic SLE$_{8/3}$ 
results and the simulations for the four random variables are shown in 
figures \ref{rv1} to \ref{rv4}. 
We plot the differences for $N=200K, 500K, 1000K$. The most important 
feature of these plots is the scale on the vertical axis. 
The full vertical scale on each of the four plots is only $3 \times 10^{-4}$. 
For the most part the deviation of the SAW simulation results from the 
exact results for SLE$_{8/3}$ comes from statistical errors, i.e., from 
the fact that we cannot run the Monte Carlo simulation forever. 
Systematic errors from the finite length of the SAW's can be seen in figures
\ref{rv2} and \ref{rv3} for values of the random variable just above $1$. 

\begin{table}
\begin{center}
\begin{tabular}{|r|r|r|r|r|}
\hline
N & samples  & $b-5/8$ & $\bp-5/48$  \\
\hline
    100K  &  566M &  -0.0012195299 &  -0.0003553937 \\
    200K  &  399M &  -0.0009095889 &  -0.0003603211 \\
    500K  &  254M &  -0.0004963860 &  -0.0002366523 \\
    1000K &  180M &  -0.0003692820 &  -0.0002336083 \\
\hline
\end{tabular}
\caption{\leftskip=25 pt \rightskip= 25 pt 
The estimates of $b$ and $\bp$ from the SAW simulations. 
We use four different lengths of walks. The second column gives 
the number of samples used in millions.
}
\label{table_exps} 
\end{center}
\end{table}

Let $W$ be one of the random variables $X,Y,R,S$ and 
$\phi_w(z)$ the corresponding conformal map where $w=x,y,r,s$. 
Taking the $\log$ of \reff{lsw_formula} we have 
\bea
\log[\probNp( W \le w)] = 
b \, \log |\phi_w^\prime(0)| + \bp \, \log |\phi_w^\prime(i)|
\label{lsfit}
\eea 
The right side is linear in $b$ and $\bp$, so we can do a least 
squares fit to estimate these two exponents. 
For each of the four random variables we take a discrete set of values 
of $w$. For $X$ these values range from $0$ to $5$ by $0.01$. 
For $Y,R,S$ they range from $1$ to $5$ by $0.01$. 
Our simulations estimate 
$\log[\probNp( W \le w)]$. The calculations in section \ref{exact_radial}
give $\log|\phi_w^\prime(0)|$ and $\log|\phi_w^\prime(i)|$. 
The results of a least square fit for $b$ and $\bp$ using these 
$1900$ cases of \reff{lsfit} are shown in table \ref{table_exps} 
for SAW's with 100K, 200K, 500K and 1000K steps. 
The values of $b$ and $\bp$ are given in terms of their difference
with the conjectured values of $b=5/8$ and $\bp=5/48$ \cite{lsw_saw}.
Note that values of $b$ and $\bp$ we obtain from the simulations are 
within a few hundredths of a percent of the conjectured values.

\clearpage

\section{Conclusions}
\label{conclusions}

We have shown how one may obtain radial SLE$_{8/3}$ from the scaling 
limit of the fixed-length SAW in the half-plane. We take a curve 
from the scaling limit and apply the conformal automorphism of 
the half-plane that fixes the origin and maps the endpoint of the 
curve to $i$. This transformation by itself does not transform 
the fixed-length SAW into radial SLE$_{8/3}$ from from $0$ to $i$. 
We must also change the probability measure by weighting each curve
from the fixed-length scaling limit by $R^p$ where $R$ is the distance of 
the endpoint of the curve to the origin. The power $p$ is conjectured
to be $(\rho-\gamma)/\nu=-61/48$.

Our heuristic derivation of the conjectured relationship 
of the fixed-length SAW to radial SLE$_{8/3}$ is further supported 
by Monte Carlo simulations of the SAW. In particular, we 
computed estimates of the scaling exponents $b$ and $\bp$ and found 
values that agree with the conjectured values within 
a few hundredths of a percent.
Our simulations also gave a partial test of predictions of the 
SLE partition function \reff{explicit_sle_pf} for the SAW in the half-plane 
with arbitrary length. 

An obvious open problem is to prove the relationship
between the fixed-length SAW and radial SLE$_{8/3}$. Given the near 
total absence of rigorous results on the two-dimensional SAW, 
progress on this problem would be a major breakthrough. 

\ifodd \verbose

\section{Code}

\subsection{Fast RV's}

Let $x+iy$ be the endpoint of $\omega$. The conformal automorphism
of $\half$ that maps $x+iy$ to $i$ is given by 
\bea
\phi(z)={yz \over y^2+x^2-xz}
\eea

The brute force approach to computing RV's is to loop through all the 
points on the walk $\omega$. To speed this up we take advantage of the 
fact that since the walk is nearest neighbor, jumping $l$ steps ahead
in the walk can only move you a distance at most $l$. 
Under the conformal map $\phi$ this distance is changed. To bound 
the distance traveled along $\phi(\omega)$ we need a bound on $\phi^\prime$. 

We have
\bea
\phi^\prime(z)={y(x^2+y^2) \over (y^2+x^2-xz)^2}
\eea
We need to bound this in a neighborhood of $z$. 
\bea
\phi^\prime(z+\delta)={y(x^2+y^2) \over (y^2+x^2-xz-x \delta)^2}
\eea
The triangle inequality gives 
\bea
|y^2+x^2-xz-x \delta| \ge |y^2+x^2-xz| - |x| |\delta|
\eea
So if we restrict $\delta$ by 
\bea
|\delta| \le {|y^2+x^2-xz| \over 2 |x|},
\label{restrict}
\eea
then 
\bea
|\phi^\prime(z+\delta)| \le 4 {|y(x^2+y^2)| \over |y^2+x^2-xz|^2} = 
4 |\phi^\prime(z)|
\eea
We have $|\phi(z)-\phi(z+\delta)| \le \max |\phi^\prime(w)| |\delta|$, where
the max is over the line segment from $z$ to $z+\delta$. So 
assuming \reff{restrict} holds, 
$|\phi(z)-\phi(z+\delta)| \le 4 |\phi^\prime(z)| |\delta|$. 
Let $d$ be the distance we can travel in the 0-i plane before we need 
to check the RV again. We want 
$|\phi(z)-\phi(z+\delta)| \le d$, so it suffices to 
have $4 |\phi^\prime(z)| |\delta| \le d$. Thus our two conditions on 
$\delta$ are 
\bea
|\delta| &\le& {|y^2+x^2-xz| \over 2 |x|} \nonumber \\
|\delta| &\le& { d \over 4 |\phi^\prime(z)|}  \nonumber \\
\eea
We can then take $l$ to be the largest integer less than $|\delta|$.

\subsection{Least squares fit}

study\_rv\_all.c produces file (temp/den\_err0) with the density and error
bars in a format that gnuplot uses. This file is input for 
study\_saw\_radial\_sle.c.  This program produces two files.
Suggested name temp/least\_sq.plt is for plotting with gnuplot with error
bars. Format is log(sin(theta)),log(prob), log(prob\_low), log(prob\_high). 
It should be approximately linear. 
Suggested name temp/least\_sq.dat is for input to simple\_least\_sq.c.
Format is log(sin(theta)),log(prob),stand dev.
simple\_least\_sq.c. will print the crucial info (estimate of the 
slope) and produces plot files.

\else
\fi

\bigskip

\noindent {\bf Acknowledgments:}
Oded Schramm's suggestion to test if simply applying the conformal map 
to the fixed-length SAW would give radial SLE$_{8/3}$ stimulated the 
author's interest in this problem. The author has also benefited
from discussions with Greg Lawler and Wendelin Werner.
This research was supported in part by the National Science Foundation 
under grant DMS -0758649.


\begin{thebibliography}{}

\bibitem{bb_review} M.~Bauer, D.~Bernard,
2D growth processes: SLE and Loewner chains, Phys. Rep.
{\bf 432}, 115-221 (2006).
Archived as arXiv:math-ph/0602049v1.

\bibitem{cardy_review}
J.~Cardy, SLE for theoretical physicists, Ann. Physics {\bf 318}, 81-118 (2005).
Archived as arXiv:cond-mat/0503313v2 [cond-mat.stat-mech].

\bibitem{clisby}
N.~Clisby,
Efficient implementation of the pivot algorithm for self-avoiding walks,
J. Statist. Phys. {\bf 140}, 349-392 (2010).
Archived as arXiv:1005.1444v1 [cond-mat.stat-mech].

\bibitem{dc_smirnov} 
H.~Duminil-Copin,S.~Smirnov,
The connective constant of the honeycomb lattice equals $\sqrt{2+\sqrt{2}}$.
Preprint, 2010. Archived arXiv:1007.0575v1 [math-ph].

\bibitem{saw_bridges} B.~Dyhr, M.~Gilbert, T.~Kennedy, G.~Lawler, S.~Passon, 
The self-avoiding walk in a strip. Preprint (2010).
Archived as arXiv:1008.4321v1 [math.PR].

\bibitem{flory}
P.J. Flory, The configuration of a real polymer chain, J. Chem. Phys. {\bf 17},
303-310 (1949). 

\bibitem{kn_review} W.~Kager, B.~Nienhuis, 	
A guide to stochastic Loewner evolution and its applications, 
J. Statist. Phys. {\bf 115}, 1149-1229 (2004).
Archived as arXiv:math-ph/0312056v3.

\bibitem{kennedy_pivot} 
T.~Kennedy,  A faster implementation of the pivot algorithm 
for self-avoiding walks,  J. Stat. Phys.  {\bf 106}, 407-429 (2002).
Archived as arXiv:cond-mat/0109308v1.

\bibitem{Kennedya} 
T.~Kennedy, 
Monte Carlo tests of SLE predictions for 2D self-avoiding walks, 
Phys. Rev. Lett. {\bf 88}, 130601 (2002).
Archived as arXiv:math/0112246v1 [math.PR]. 

\bibitem{Kennedyb} 
T.~Kennedy,	
Conformal invariance and stochastic Loewner evolution predictions 
for the 2D self-avoiding walk - Monte Carlo tests,
J. Stat. Phys. {\bf 114}, 51--78 (2004).
Archived as arXiv:math/0207231v2 [math.PR].


\bibitem{kestena} 
H.~Kesten, On the number of self-avoiding walks,  
J. Math. Phys  {\bf 4}, 960-969  (1963). 

\bibitem{kestenb} 
H.~Kesten, On the number of self-avoiding walks II,  
J. Math. Phys  {\bf 5}, 1128-1137 (1964). 

\bibitem{lawler_part_func} 
G. Lawler,  Partition functions, loop measure, and versions of SLE.
J. Statist. Phys. {\bf 134}, 813-837 (2009).

\bibitem{lawler_book} 
G. Lawler,  {\it Conformally Invariant Processes in the Plane}.  
American Mathematical Society (2005).

\bibitem{lawler_utah} 
G.~Lawler, Schramm-Loewner evolution, in {\it Statistical Mechanics}, 
S. Sheffield and T. Spencer, ed., IAS/Park City Mathematical Series, 
AMS, 231-295 (2009). 
Archived as arXiv:0712.3256v1 [math.PR].

\bibitem{lsw_saw} 
G.~Lawler, O.~Schramm, and W.~Werner, 
On the scaling limit of planar self-avoiding walk, 
{\it Fractal Geometry and Applications: a Jubilee of Benoit Mandelbrot, 
Part 2}, 339--364, {\em Proc. Sympos. Pure Math. 72}, 
Amer. Math. Soc., Providence, RI, 2004.
Archived as arXiv:math/0204277v2 [math.PR].

\bibitem{lsw_restriction} 
G.~Lawler, O.~Schramm, W.Werner, 
Conformal restriction: the chordal case,
J. Amer. Math. Soc. {\bf 16}, 917­-955 (2003).
Archived as arXiv:math/0209343v2 [math.PR].

\bibitem{ms_book}  
N. Madras and G. Slade, {\em The Self-Avoiding Walk}.  
Birkh\"{a}user (1996). 

\bibitem{nienhuis}
B.~Nienhuis, Exact critical exponents for the $O(n)$ models in two dimensions,
Phys. Rev. Lett. {\bf 49} 1062-1065 (1982).

\bibitem{schramm} 
O.~Schramm, Scaling limits of loop-erased random walks and 
uniform spanning trees, Israel J. Math. {\bf 118}, 221-288 (2000).
arXiv:math/9904022v2 [math.PR].

\bibitem{schramm_private} O.~Schramm, private communication, November, 2002.

\bibitem{werner_review} 
W.~Werner, Random planar curves and Schramm-Loewner evolutions, 
Ecole d'Et\'e de Probabilit\'es de Saint-Flour XXXII - 2002, 
Lecture Notes in Mathematics 1840, Springer-Verlag, 107–195 (2004).
Archived as arXiv:math/0303354v1 [math.PR].

\end{thebibliography}
\end{document}